\documentclass[12pt]{amsart}
\usepackage{amscd}
\usepackage{amssymb}
\newcommand{\C}{\mathbb C}
\newcommand{\Q}{\mathbb Q}
\newcommand{\Z}{\mathbb Z}
\newcommand{\X}{\mathcal X}
\usepackage[curve]{xypic}
\newcommand{\comment}[1]{\marginpar{\sffamily{\noindent\tiny #1
   \par}\normalfont}}
\renewcommand{\comment}[1]{}
\newbox\mybox
\def\overtag#1#2#3{\setbox\mybox\hbox{$#1$}\hbox to
  0pt{\vbox to 0pt{\vglue-#3\vglue-\ht\mybox\hbox to \wd\mybox
      {\hss$\scriptstyle#2$\hss}\vss}\hss}\box\mybox}
\def\undertag#1#2#3{\setbox\mybox\hbox{$#1$}\hbox to 0pt{\vbox to
    0pt{\vglue#3\vglue\ht\mybox\hbox to \wd\mybox
      {\hss$\scriptstyle#2$\hss}\vss}\hss}\box\mybox}
\def\lefttag#1#2#3{\hbox to 0pt{\vbox to 0pt{\vss\hbox to
      0pt{\hss$\scriptstyle#2$\hskip#3}\vss}}#1}
\def\righttag#1#2#3{\hbox to 0pt{\vbox to 0pt{\vss\hbox to
      0pt{\hskip#3$\scriptstyle#2$\hss}\vss}}#1}

\def\Dot{\lower.2pc\hbox to 2.5pt{\hss$\bullet$\hss}}
\def\Circ{\lower.2pc\hbox to 2.5pt{\hss$\circ$\hss}}
\def\Vdots{\raise5pt\hbox{$\vdots$}}
\def\splicediag#1#2{\xymatrix@R=#1pt@C=#2pt@M=0pt@W=0pt@H=0pt}

\renewcommand\frame[2][3pt]{\hbox{$\vcenter{\hbox{\vrule\vbox 
{\hrule\kern#1\hbox{\kern#1$#2$\kern#1}\kern#1\hrule}\vrule}}$}}
\newcommand\lineto{\ar@{-}}
\newcommand\dashto{\ar@{--}}
\newcommand\dotto{\ar@{.}}
\newcommand{\bt}{\bullet}
\newtheorem{theorem}{Theorem}[section]

\newtheorem*{theorem*}{Theorem}
\newtheorem{proposition}[theorem]{Proposition}
\newtheorem{lemma}[theorem]{Lemma}
\newtheorem{corollary}[theorem]{Corollary}
\newtheorem*{corollary*}{Corollary}

\newtheorem*{conjecture*}{Conjecture}
\newtheorem*{BStheorem*}{Bhupal-Stipsicz Theorem}

\theoremstyle{definition}
\newtheorem{remark}[theorem]{Remark}
\newtheorem*{example*}{Example}
\newtheorem*{examples*}{Examples}
\newtheorem*{remarks*}{Remarks}
\newtheorem{definition}[theorem]{Definition}

\begin{document}
\title[Log-terminal smoothings of graded normal surface singularities]
{Log-terminal smoothings of graded normal surface singularities}
\author{Jonathan Wahl}
\thanks{Research supported under NSA grant no.\ H98230-10-1-0365}
\address{Department of Mathematics\\The University of North
  Carolina\\Chapel Hill, NC 27599-3250} \email{jmwahl@email.unc.edu}
\begin{abstract} Recent work (\cite{SSW}, \cite{bs}) has produced a complete list of weighted homogeneous surface singularities admitting smoothings whose Milnor fibre has only trivial rational homology (a ``rational homology disk").  Though these special singularities form an unfamiliar class and are rarely even log-canonical, we prove the
\begin{theorem*} A rational homology disk smoothing of a weighted homogeneous surface singularity can always be chosen so that the total space is log-terminal.  In particular, this smoothing is $\Q$-Gorenstein.
\end{theorem*}
The key idea is to define a finite ``graded discrepancy" of a normal graded domain with $\Q$-Cartier canonical divisor, and to study its behavior for a smoothing.  

\end{abstract}
\maketitle
    
\begin{center}{\textbf{Introduction}}
    \end{center}
    
    Let  $(X,0)$ be the germ of an isolated complex normal singularity.  Suppose that its canonical module $K_X$ is $\Q$-Cartier; thus, the canonical sheaf $K_{X-\{0\}}$ has some finite order $m$.  The  \emph{index one} (or \emph{canonical}) cover $(T,0)\rightarrow (X,0)$ is obtained by normalizing the corresponding cyclic cover.  $(T,0)$ has an isolated normal singularity with $K_T$ Cartier.  If $T$ is Cohen-Macaulay, then it is Gorenstein; we call such an $(X,0)$ \emph{ $\Q$-Gorenstein}.  (Warning: some authors require only that $X$ be Cohen-Macaulay with $K_X$ $\Q$-Cartier.)  That $T$ be Cohen-Macaulay is automatic if $X$ has dimension $2$, but not in general, even if $X$ itself is Cohen-Macaulay.  In fact, A. Singh (\cite{si}, 6.1) constructs an example of an isolated $3$-dimensional \emph{rational} (hence Cohen-Macaulay) singularity $X$ with $K_X$ $\Q$-Cartier whose index $1$ cover is not Cohen-Macaulay.
    
    Now suppose $(X,0)$ is a $\Q$-Gorenstein normal surface singularity (e.g., a rational singularity).     We will say that a smoothing $f:(\mathcal X,0)\rightarrow (\C,0)$ of $X$ is \emph{$\Q$-Gorenstein} if it is the quotient of a smoothing of the index one cover of $X$.  The basic examples are certain smoothings of the cyclic quotient singularities of type $rn^2/rnq-1$, first mentioned in \cite{lw}, (5.9).   In this particular case, $\mathcal X$ is a cyclic quotient of $\C^3$---it is even a \emph{terminal} singularity--- and a cyclic quotient by a smaller group is the total space of the smoothing of the index one cover of $X$ (which is an $A_{rn-1}$-singularity).  The importance of $\Q$-Gorenstein smoothings was first noticed in the work of Shepherd-Barron and Koll\'ar \cite{ksb}.
    But even if the (now three-dimensional) total space $\mathcal X$ of a smoothing of $X$ has $K_{\mathcal X}$ $\Q$-Cartier, it does not immediately follow that the smoothing is $\Q$-Gorenstein.
           
       The main point of this work is that in an important special case, one can deduce that a smoothing is $\Q$-Gorenstein by proving the much stronger result that the total space $\mathcal X$ can be chosen to be \emph{log-terminal}.   The result is surprising since the original singularities $X$ are generally not even log-canonical, i.e., have discrepancies $-\infty.$   
       
       In the early 1980's the author constructed many examples (both published and unpublished) of surface singularities that admit smoothings with Milnor number $0$, i.e., whose Milnor fibre is a $\Q$HD (rational homology disk).  The simplest of these were the aforementioned $n^2/nq-1$ cyclic quotients.  These smoothing examples were constructed in one of two ways.   The first was an explicit ``quotient construction" (\cite{w3}, 5.9).  One starts with (a germ of) an isolated $3$-dimensional Gorenstein singularity $(\mathcal Z,0)$; a finite group $G$ acting on it, freely off the origin; and a $G$-invariant function $f$ on $\mathcal Z$, whose zero locus $(W,0)$ has an isolated singularity (hence is normal and Gorenstein).  Then  $f:(\mathcal Z,0)\rightarrow (\mathbb C,0)$ is a smoothing of $W$, with a Milnor fibre $M$, which in these example is simply-connected.   One also has  $f:(\mathcal Z/G,0)\rightarrow (\mathbb C,0)$,  giving a $\mathbb Q$-Gorenstein smoothing of $(W/G,0)\equiv(X,0)$, with Milnor fibre $M/G$ (note $G$ acts freely on $M$).   Here the order of $G$ equals the Euler characteristic of $M$, so $M/G$ has Euler characteristic $1$, hence is a $\mathbb Q$HD Milnor fibre.   In all these examples, $\mathcal Z$ is actually canonical Gorenstein, hence the total space $\mathcal Z/G$ of the smoothing is log-terminal.
       
       But, other examples could be constructed only using using H. Pinkham's method of ``smoothing with negative weight"  for a weighted homogeneous singularity \cite{p3}.  For these cases (e.g., the triply-infinite family of type $\mathcal M$ of \cite{SSW}, (8.3)), the properties of the total space of the smoothing are much less obvious, and it is not at all clear that the smoothings are $\Q$-Gorenstein.  
       
       The problem of finding three-manifolds which ``nicely bound" rational homology disks became of interest in symplectic topology, and led to a proof that our old list of examples was complete in the weighted homogeneous case (\cite{SSW} and \cite{bs}).   Thus, one knows the resolution dual graph of all weighted homogeneous surface singularities admitting a $\Q$HD smoothing.   Nonetheless, this explicit class of singularities is rather mysterious.  We show that such smoothings are $\Q$-Gorenstein, and can even be assumed to have log-terminal total space.  In (3.4), we prove 
        %The resolution dual graph of $W/G$ can be calculated by resolving $\mathcal Z$ or $W$ and dividing by the induced action of $G$.  
  \begin{theorem*}  Let $(X,0)$ be a weighted homogeneous surface singularity admitting a $\Q$HD smoothing.  Then
 \begin{enumerate} 
 \item A $\Q$HD smoothing occurs over a one-dimensional smoothing component of $(X,0)$.
  \item Let $f:(\mathcal X,0)\rightarrow (\C,0)$ be the induced smoothing over the normalization of that component.  Then $(\mathcal X,0)$is log-terminal, and the smoothing is $\Q$-Gorenstein.
 \end{enumerate}
  \end{theorem*} 
  
  The first statement is  \cite{w5}, Corollary 8.2 (one guesses that the smoothing component is smooth, so normalization is unnecessary); the main point is the second assertion.  We do \emph{not} claim that log-terminality remains true after base-change (see  Remark (3.5)(1) below). 
     
       Our main technical tool is to define and study a \emph{graded discrepancy} $\alpha(X)$ for an isolated normal graded singularity $X$ with $K_X$  $\Q$-Cartier.   The blow-up of the weight filtration $Z\rightarrow X$ has irreducible exceptional fibre $E$, so that $(Z,E)$ is a log-terminal pair.  Then we define ((1.2) below) $\alpha(X)$ to be the discrepancy of $K$ along $E$.  This invariant depends not only on $X$ but upon the grading as well.  A key observation (Proposition 1.4) is:   $$ \alpha(X)>-1\  \text{implies}\ X\ \text{is\ log-terminal}. $$ 
       
       A basic result relates the graded discrepancy of a singularity to that of a graded hypersurface section (see (1.5) and (1.6)):
       
       \begin{theorem*} Let $\X=\text{Spec}\ A$ be a graded normal isolated singularity, with $K_{\X}$ $\Q$-Cartier, whose minimal set of generators $z_1,z_2,\cdots,z_t$ has weights $m_1,m_2,\cdots,m_t$, with GCD equal $1$.  %, and GCD$(m_i,m_t)=1$ for all $i<t$.  
    Suppose  $X=\text{Spec}\ A/(f)$ is normal and graded, with isolated singularity. % Assume that GCD$(m_1,\cdots,m_{t-1})=1$.  
    Assume further one of the following:
    \begin{enumerate}
    \item $f=z_t$, and $m_1,\cdots,m_{t-1}$ have GCD equal $1$.
    \item $f\in m_A^2$.
    \end{enumerate}
    Then $\alpha(\X)=\alpha(X)+\text{weight}\ f.$
  \end{theorem*}     
 \begin{corollary*} Under the above hypothesis, if in addition $X$ is a $\Q$-Gorenstein surface singularity with $\alpha(X)>-2$, then $\X$ is log-terminal (and in particular the corresponding smoothing is $\Q$-Gorenstein).
 \end{corollary*}
 If $X$ is a weighted homogeneous surface singularity, we show in Section $2$ how to compute $\alpha(X)$ in terms of some topological invariants $\chi$ and $e$, introduced by W. Neumann.   Proposition 2.1 shows $$\alpha(X)=-1-(\chi/e).$$ But \cite{SSW} limits greatly the possible resolution dual graphs of the singularities with a $\Q$HD smoothing; computing $\chi$ and $e$, we conclude $\alpha(X)>-2$ in those cases.  (While \cite{bs}  pinned down the exact list of graphs, that more precise result is not needed in this proof).
Thus $\X$ is the total space of a smoothing of $X$ and $\alpha(\X)>-1$ by the last Theorem, so it is log-terminal.  
         
% Because of the relation of the graded discrepancies, we can show $\mathbb X$ is log-terminal (i.e., $\alpha(\mathbb X)>-1$) once we show  that the $X$ in the theorem have $\alpha(X)>-2.$  But we can express $\alpha(X)$ for a graded normal surface singularity  
  Finally, we apply an old result of K. Watanabe \cite{wat2} in Corollary 4.6, showing exactly which weighted homogeneous surface singularities are $\Q$-Gorenstein (an analytic condition).  The invariant $\chi/e$ shows up naturally. 

We have profited from conversations with J\'anos Koll\'ar and S\'andor Kov\'acs.

\bigskip

     \section{Seifert partial resolution and graded discrepancy}
     \bigskip
     Let $z_{1},\dots,z_t$ be coordinate functions on an affine space $\C
^{t}$, where $z_{i}$ has positive integer weight $m_{i}$.  Assume that GCD$(m_1,\cdots,m_t)=1$.
Blowing-up
the corresponding weight filtration gives the \emph{weighted} or \emph{filtered blow-up}
$\pi:\mathcal Z\rightarrow \C^t$, an isomorphism off $\pi^{-1}(0)\equiv \mathcal E$, an irreducible Weil divisor isomorphic to the corresponding weighted projective space.   $\mathcal Z$ is covered by $t$ affine varieties $U_{i}$,
each of
which is a quotient of an affine space $V_{i}$ by a cyclic group of order
$m_{i}$.  Consider on $V_{1}$  coordinates $x,y_2,\cdots,y_t$,
 related to the $z_{i}$ via
$$z_{1}=x^{m_{1}},~ z_{2}=x^{m_{2}}y_{2},~ \dots,~
z_{t}=x^{m_{t}}y_{t}.$$ 
$U_{1}$ is the quotient of  $V_{1}$ by  the action of the
cyclic
group generated by $$S=[-1/m_{1},m_{2}/m_{1},\dots,m_{t}/m_{1}]=(1/m_1)[-1, m_2,\cdots,m_t]\ ,$$
where we are using the notation
$$[q_1,\dots,q_t]:=(\exp(2\pi i q_1),\dots,\exp(2\pi iq_t))\,.$$
The hyperplane $\{x=0\}$ in $V_1$ maps onto the exceptional fibre $\mathcal E\cap U_1$.

  The group $\langle S \rangle$ acts on $V_1$ without pseudo-reflections.  We describe the locus where the action is not free.  For each prime $p$ dividing $m_1$, let $J_p=\{j|\ 2\leq j \leq t,\  (p,m_j)=1\}$; it is non-empty.  If $j\notin J_p$, then $S^{m_1/p}$ acts trivially in the $j$th slot.  Thus $S^{m_1/p}$ acts trivially on the linear subspace $V_{1,p}$ defined by the vanishing of $x$ and the $y_j$ with $j\in J_p$.  The codimension of $V_{1,p}$ in $V_1$ is $1+\#J_p$, and the union of the images in $U_1$ gives the singular locus there.    
 
    Let $X=\text{Spec} \ A$ be an isolated normal singularity with good $\C^*$-action.  Then the graded domain $A=\oplus A_k$ can be written  as the quotient of a graded polynomial ring $\C[z_1,\cdots,z_t]$,  where the weights $m_i$ of the $z_i$ are assumed to have GCD $1$.  Assume further that one has chosen a minimal set of generators, i.e., that the embedding dimension is $t$.   Denote the weight filtration by $I_s=\oplus _{k\geq s}A_k$, and consider $$\oplus_{s=0}^{\infty}\ I_s u^s\subset A[u].$$   The blow-up of the weight filtration of $A$ is $$\pi:Z=\text{Proj}  \oplus_{s=0}^{\infty}I_s u^s \rightarrow X=\text{Spec} \ A,$$ and is sometimes called the \emph{Seifert partial resolution} of $X$ (of course, it depends on the choice of grading).  It is an isomorphism off $\pi^{-1}(0)\equiv E$, an irreducible Weil divisor isomorphic to $\text{Proj}\ A$.  
    
    Further, $\pi:Z\rightarrow X$ is equal to the proper transform in $\mathcal Z$ of the subvariety $\text{Spec}\ A\subset \C^t$, so $Z$ is covered by closed subvarieties $\bar{U_i}$ of the affines $U_i$ above.   More precisely, let $\{g_\alpha(z_1,\cdots,z_t)\}$ be homogeneous generators of the graded ideal defining $X$.   Then the proper transform $\bar{V}_1$ of $X$ in $V_1$ is the affine subvariety with coordinate ring $\C[x,y_2,\cdots,y_t]/(g_\alpha(1,y_2,\cdots, y_t))\equiv\mathcal A_1$, which is  a polynomial ring in $x$ over the subring generated by $y_2,\cdots ,y_t.$  Note that $\bar{V}_1$ is smooth, because $X$ has an isolated singularity (for this and other details, consult e.g. Flenner \cite{fl} or \cite{w4},(2.1)); thus, $(Z,E)$ is a log-terminal pair, locally the cyclic quotient of a smooth space plus normal crossings divisors .  In particular, the corresponding affine open subset of $Z$ is $\bar{U}_1=\bar{V}_1/\langle S \rangle$.   Since the $z_i$ form a minimal set of generators, all $y_i$ are non-$0$ in $\mathcal A_1$.  So the group acts on $\bar{V}_1$ ``without pseudo-reflections" , i.e. acts freely off the intersections $\bar{V}_1\cap V_{1,p}\equiv \bar{V}_{1,p}$, necessarily of codimension at least $2$ (but possibly less than $1+\#J_p$).
    
    We wish to compare the Seifert partial resolution $(Z,E)\rightarrow (X,0)$ with that of an appropriate hypersurface section.  So, let $f\in A_d$ be a homogeneous element with the property that the hypersurface $X'=\text{Spec}\ A/(f)\subset X$ is also normal, with isolated singularity.  Then $f=0$ on $Z$ consists of the proper transform $Z'$ of $X'$, plus the set $E'=E\cap Z'$ (which is irreducible, isomorphic to $\text{Proj}\ A/(f)$).  We are interested in achieving the condition
    $$(*)\ \ E' \  \text{is not contained in the singular locus of Z}.$$
    When $(*)$ is satisfied, one has $$(**)\ \text{generically}\   E'\  \text{is a Cartier divisor on both}\ E\  \text{and}\ Z'.$$   To verify $(*)$, it suffices to check in some $\bar{U}_i$ for which the pull-back of $E'$ is non-empty.  The singular locus is the image of the fixed-point locus in $\bar{V}_i$.  So, in say $\bar{V}_1$,  write $f=x^df(1,y_2,\cdots,y_t)\equiv x^d\bar{f}$.  Assume that $(\bar{f},x)$ is not the unit ideal.  The condition $(*)$ means that the set  $x=\bar{f}=0$ is not contained in any $\bar{V}_{1,p}$.
    
    \begin{example*} Let $X=\C^3$, with coordinates $z_1,z_2,z_3$ having weights $2,2,1$ respectively.  Then $V_1$ has coordinates $x,y_2,y_3$, and the action of $S$ is $(1/2)[-1,2,1]$; so $U_1$ has coordinates $u=x^2, \ v=xy_3,\ w=y_3^2, \ y_2$.  In other words, $U_1$ is a line cross the ordinary double point given by $uw=v^2$, with exceptional fibre defined by $u=v=0$, and singular locus the line cross the point $u=v=w=0.$  Consider now the proper transform of $f=z_3$ in $U_1$, defined by $v=w=0$.  Its intersection with the exceptional fibre is exactly the singular locus of $U_1$.  Thus, $(*)$ is not satisfied.
    %Let $X=\text{Spec}\ \C[z_1,z_2,z_3,z_4]/(z_1^3+z_2^3+z_3^3+z_4^{3m})$, where $m\geq 1$.  The weights are $m,m,m,1$.  In the previous notation, $\bar{V}_1$ is the affine variety with coordinate ring $\C[A_1,A_2,A_3,A_4]/(1+A_2^3+A_3^3+A_4^{3m})$, and the group action is given by $S=(1/m)[-1,0,0,1]$.  Suppose $m>1$.  Then the fixed locus of  any power of $S$ is given by $A_1=A_4=0$.  The function $f=z_4$ on $X$ cuts out an isolated normal singularity, but it factors $A_1\cdot A_4=A_1\bar{f}$ on $\bar{V}_1$.  While $(A_1,A_4)$ is not the unit ideal, it contains (and in fact is equal to) the ideal of  the fixed locus.  Thus, $(*)$ is not satisfied when $m>1$.
    \end{example*}
    
    The problem in the last example was that the weights other than that of $z_3$ were not relatively prime.  One does have the
    
    \begin{proposition} \label{pr:wts} Let $X=\text{Spec}\ A$ be a graded normal isolated singularity, with the minimal set of weights $m_1,m_2,\cdots,m_t$ with GCD equal $1$.  %, and GCD$(m_i,m_t)=1$ for all $i<t$.  
    Suppose  $X'=\text{Spec}\ A/(f)$ is normal and graded, with isolated singularity. % Assume that GCD$(m_1,\cdots,m_{t-1})=1$.  
    Then Condition $(*)$ is satisfied in each of the following cases:
    \begin{enumerate}
    \item $f=z_t$ and $m_1,\cdots,m_{t-1}$ have GCD equal $1$.
    \item $f\in m_A^2$.
    \end{enumerate}
    \begin{proof}  Consider the first case.  As indicated above, it suffices to look in the $\bar{V}_i$.  For $i=t$, the corresponding $\bar{f}$ is equal to $1$, so $\bar{U}_t\cap E'=\emptyset.$  So, one may restrict to, say, $\bar{V}_1$.   Here, $\bar{f}=y_t.$ 
    
    Now, $A/(z_t)$ is graded normal with isolated singularity, \emph{and} its minimal set of generators has relatively prime weights are $m_1,\cdots,m_{t-1}$.  Thus, the relevant  affine $V_1'$ has coordinate ring $\mathcal A_1/(y_t)$, with an induced action by $\langle S \rangle$ which is still free off a set of codimension at least $2$.
    
    But suppose there exists a prime $p$ so that the locus $x=y_t=0$ is contained in $\bar{V}_{1,p}$.  Then $S^{m/p}$ would act trivially on a divisor of $V_1'$.  This is a contradiction.
    
    The second case is proved similarly.
    
    \end{proof}
    
    \end{proposition}

     Now suppose $K_X$ is $\Q$-Cartier, so that for some positive integer $n$ one has on $U=X-\{0\}$ that $\omega_U^{\otimes n}\cong \mathcal O_U$.   Then $nK_Z$ has some integral coefficient $c$ along $E$, and we make the the obvious
     
     \begin{definition} \emph{The graded discrepancy $\alpha(X)$ is $c/n$}.  
     \end{definition}
     
    \begin{remark} Note that $\alpha(X)$ depends upon the grading.  For instance, it is easy to check that  if $\C^2$ has  coordinates with weights $1$ and $m$, then $\alpha =m.$  
    \end{remark}
    
    However, one has the
     
     \begin{proposition}  If $\alpha(X)>-1$, then $X$ is log-terminal.  If $\alpha(X)\geq -1$, then $X$ is log-canonical.
     \begin{proof}  The key point is that the pair $(Z,E)$ is log-terminal.  So, if $\alpha(X)\geq -1$, resolving singularities will not give any exceptional components whose discrepancy is less than $\alpha(X)$.  (See \cite{k-m}, Section $5.2$ for this and other facts about log-terminal singularities, such as that they are cyclic quotients of rational Gorenstein singularities, hence Cohen-Macaulay.)
     \end{proof}
     \end{proposition}
     
     On the other hand, even when $\alpha(X)$ is smaller than $-1$, it can be a useful invariant even though the usual discrepancy of the singularity is $-\infty$.
     
     \begin{example*} Let $Y$ be a smooth projective variety, $L$ an ample line bundle.  Then the cone $X=\text{Spec} \oplus _{i\geq 0} \Gamma(Y,L^{\otimes i})$ is a normal graded ring with isolated singularity.  The Seifert partial resolution has total space the geometric line bundle $Z=V(L^{-1})\rightarrow Y$, and the exceptional divisor $E$ is the zero-section.   
     If  there exist integers $m$ and $n\neq 0$ so that $K_Y^{\otimes n}\cong L^{\otimes m}$, then $K_X$ is $\Q$-Cartier.  Writing $K_Z\equiv cE$ and using the adjunction formula, one easily finds that $$\alpha(X)=-1-(m/n).$$
      \end{example*}
      
      An important use of the invariant $\alpha(X)$ arises when taking graded hypersurface sections, for instance when one views $X$ as the total space of a graded smoothing of a singularity of dimension one less.
      
      \begin{theorem}  Let $\pi:Z\rightarrow X=\text{Spec}\ A$ be the Seifert partial resolution of an isolated normal singularity with good $\C^*$-action, for which $K_X$ is $\Q$-Cartier.  Suppose $f\in A_d$ is such that the Cartier divisor $X'$ defined by $f$ also has an isolated normal singularity, and satisfies Condition $(*)$.  Then
      \begin{enumerate}
      \item The Cartier divisor defined by $f$ on $Z$ is generically $dE+Z'$, where $Z'\rightarrow X'$ is the Seifert partial resolution of $X'$.
      \item The graded discrepancies are related by $$\alpha(X)=\alpha(X')+d.$$
      \end{enumerate}
      \begin{proof}  The first assertion follows from the discussion above relating weighted blow-up of a singularity to the proper transform in the weighted blowing-up of a graded affine space.  Condition $(*)$ means that generically along $E'=E\cap Z'$, $E$ and $Z$ are Cartier divisors.
      For the second assertion, one can write generically $K_Z\equiv \alpha E$.  By adjunction, generically along $E'$ one has $K_{Z'}\cong \mathcal O_{Z'}(\alpha E+Z').$  But the divisor of $f$ is trivial on $Z$, so generically $Z'\equiv -dE$.  The result follows.
      \end{proof}
      \end{theorem}
      \begin{corollary} Let $X=\text{Spec}\ A$ be an isolated normal singularity with good $\C^*$-action,  with $K_X$ $\Q$-Cartier.  Suppose  $f\in A_d$ also defines an isolated normal singularity $X'$, and $f$ satisfies either of the conditions of Proposition 1.1.  If $a(X')>-2$ (respectively, $a(X')\geq -2$), then $X$ is log-terminal (resp. log-canonical).
      \end{corollary}

      Suppose $X'$ is a singularity with good $\C^*$-action. A  deformation of $X'$ is said to have \emph {negative weight} if roughly speaking one perturbs the defining equations by terms of smaller degree.  Specifically, a smoothing of negative weight consists of an isolated singularity $X$ with good $\C^*$-action, a function $f:X\rightarrow \C$ of some positive weight $d$, and a \emph{graded} isomorphism of the special fibre with $X'$.  Even if $K_{X'}$ is $\Q$-Cartier, $K_{X}$ need not be (see Remark 3.2 below.)  However, we have the
      
      \begin{corollary} \label{co:f} Let $X'$ be an isolated normal singularity with good $\C^*$-action, and $K_{X'}$ $\Q$-Cartier.  Suppose $f:X\rightarrow \C$ gives a smoothing of $X'$ of negative weight, so that $K_X$ is $\Q$-Cartier.  If $a(X')>-2$, then \begin{enumerate}
      \item $X$ is log-terminal, in particular rational and $\Q$-Gorenstein
      \item $X'$ is Cohen-Macaulay and in fact $\Q$-Gorenstein.
      \item The smoothing is $\Q$-Gorenstein, in particular a quotient of a corresponding smoothing of the index $1$ cover of $X'$.     
         \end{enumerate}
      \begin{proof}  If $X$ and $X'$ have the same embedding dimension, then $f$ satisfies Condition (2) of Proposition 1.1.   Otherwise, $f$ is among a minimal set of generators for the graded ring of $X$, while the other generators have weights which are those of $X'$; thus, the other weights have GCD equal to $1$.  The first assertion is therefore just the preceding Corollary, combined with familiar facts about log-terminal singularities.  The index $1$ cover $T$ of $X$ has an isolated canonical Gorenstein singularity (Corollary (5.2.1) of \cite{k-m}), inducing a cover $T'\rightarrow X'$, which agrees with the index $1$ cover of $X'$ off the singular point.  But $T'$ is Cohen-Macaulay (since $T$ is), so it is normal, and thus equal to the index $1$ cover of $X'$.  
      \end{proof}
      \end{corollary}
   \bigskip
      
      \section{K for a normal graded surface singularity}
      \bigskip
We recall  briefly the topological description of normal surface singularities $X=\text{Spec}\ A$ with good $\C^*$-action; the analytic classification will be described later. 

 One has the cyclic quotient singularities of type $n/q$, where $0<q<n$, $(q,n)=1$, defined as the quotient of $\C^2$ by the cyclic group action generated by $$[1/n,q/n].$$
 Writing the continued fraction expansion $n/q=b_1-1/b_2-\cdots -1/b_s,$ each $b_i\geq 2$, the minimal resolution of the singularity consists of a string of rational curves, with dual graph
$$
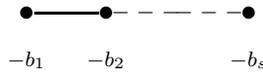
\xymatrix@R=6pt@C=24pt@M=0pt@W=0pt@H=0pt{
\\&\undertag{\bt}{-b_1}{6pt}\lineto[r]&\undertag{\bt}{-b_2}{6pt}
\dashto[r]&\dashto[r]
&\undertag{\bt}{-b_s}{6pt}\\
&&&&\\
&&&&\\
&&&&\\
}
$$

 This singularity is isomorphic to $n/q'$, where $qq'\equiv 1\  \text{mod} \ n,$ and the continued fraction expansion $n/q'$ gives the $b_i$'s in the opposite order.
 
 Every other normal surface singularity has a unique $\C^*$-action, with minimal good resolution graph $\Gamma$ which is star-shaped:
$$
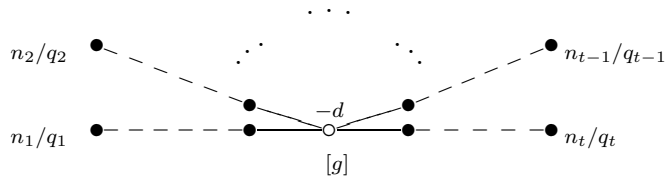
\xymatrix@R=4pt@C=24pt@M=0pt@W=0pt@H=0pt{\\
\lefttag{\bullet}{n_2/q_2}{8pt}\dashto[ddrr]&
&\hbox to 0pt{\hss\lower 4pt\hbox{.}.\,\raise3pt\hbox{.}\hss}
&\hbox to 0pt{\hss\raise15pt
\hbox{.}\,\,\raise15.7pt\hbox{.}\,\,\raise15pt\hbox{.}\hss}
&\hbox to 0pt{\hss\raise 3pt\hbox{.}\,.\lower4pt\hbox{.}\hss}
&&\righttag{\bullet}{n_{t-1}/q_{t-1}}{8pt}\dashto[ddll]\\
\\
&&\bullet\lineto[dr]&&\bullet\lineto[dl]\\
\lefttag{\bullet}{n_1/q_1}{8pt}\dashto[rr]&&
\bullet\lineto[r]&\overtag{\circ}{-d}{8pt}\undertag{}{[g]}{6pt}\lineto[r]&\bullet
\dashto[rr]&&\righttag{\bullet}{n_{t}/q_{t}}{8pt}\\&~\\&~\\&~\\&~}
$$
The strings of $\Gamma$ are described uniquely by the continued
fractions shown, starting from the node.  The central curve $C$ has genus $g$.  If $t=0$ one has a cone.  If  $g=0$, minimality and the exclusion of cyclic quotients require $t\geq 3$.  The Seifert partial resolution of the singularity contains the one exceptional curve $C$, with $t$ cyclic quotient singularities on it.

The above graph  completely determines the topology of the singularity (i.e., gives the Seifert invariants).  One has two topological invariants introduced by Neumann (e.g., \cite{ne}, p. 250):
$$e=d-\sum_{i=1}^t(q_i/n_i) $$
$$\chi=2g-2+t-\sum_{i=1}^t(1/n_i)=2g-2+\sum_{i=1}^t(1-1/n_i).$$
While $e>0$, one can easily show that $\chi<0$ exactly for quotient singularities (which are log-terminal), and $\chi=0$ exactly for the familiar list of graded log-canonical singularities.

%\begin{remark} We will see in section x.x below the role of $\chi/e$ in writing the canonical line bundle on the minimal good resolution as numerically equivalent to a rational combination of exceptional curves. 
%\end{remark}

%\section{K on the minimal good resolution}
Let $(\tilde{X},E) \rightarrow (X,0)$ be the minimal good resolution, with resolution graph as above.  The line bundle $K_{\tilde{X}}$ is numerically equivalent to a $\Q$-combination of exceptional curves, whose coefficients we examine.  For a cyclic quotient, this is an exercise using Cramer's rule (cf. \cite{lw}, 5.9.(iii)); but we shall not use this.  Exclude the cyclic quotients, and consider the $\Q$-cycle $Z\equiv-(K+E)$, writing in terms of cycles supported on the $t$ strings, plus a multiple of the central curve $C$:
$$Z=\sum_{k=1}^t  Y_k +\beta C.$$  

We solve for each $Y_k$ in terms of $n_k,\ q_k$, and $\beta$, and then solve for the appropriate $\beta.$   For any irreducible exceptional curve $F$, denote $d_F=-F\cdot F$ and $u_F=$ number of neighbors of $F$.  Thus, we have $Z\cdot F=-(d_F+2g_F-2+u_F-d_F)=2-2g_F-u_F,$ so that $Z\cdot F$ equals $1$ on the $t$ ends, $2-2g-t$ at $F=C$, and is $0$ elsewhere.

Consider first one string, of type $n/q$.   Suppose the exceptional curves counting out from $C$ are $E_1,E_2,\cdots,E_s$.  Let $e_i$ be the $\mathbb Q$-cycle supported on the string defined by $e_i(E_j)=-\delta_{ij}$.  Then the cycle $Y=\beta e_1-e_s$ dots to $-\beta$ with $E_1$, to $1$ with $E_s$, and $0$ with the intermediate $E_j$.   It is easy to see that the $\mathbb Q$-cycle $e_1$ has coefficient $1/n$ at $E_s$, and its coefficients increase strictly as one goes to $E_1$, where it is $q/n$.  By the same argument, $e_s$ has coefficient $1/n$ at $E_1$ , increasing to $q'/n$ at $E_s$, where as usual $qq'\equiv 1$ mod $n$.  In particular, the coefficient of $Y$ at $E_1$ is $(\beta q -1)/n,$ and the coefficient at the end $E_s$ is $(\beta-q')/n.$  Note that if $s=1$, then $q=q'=1$, so the formulas are still correct.

Choosing the corresponding $Y_k$ for each string, we show that for appropriate $\beta$, the cycle $Z=\sum_{k=1}^t  Y_k +\beta C$ represents $-(K+E)$.  Given the coefficient of each $Y_k$ at the curve intersecting $C$, the condition $Z\cdot C=2-2g-t$ becomes $$-\beta d+\sum_{k=1}^t (\beta q_k-1)/n_k =2-2g-t.$$  We therefore must set $$\beta=\chi/e.$$

\begin{proposition}  Writing $K_{\tilde{X}}$ numerically as a rational combination of exceptional curves, the coefficient of the central curve $C$ is $-1-(\chi/e)$.
\end{proposition}
%\begin{remark}  It is a simple exercise to list all the graphs with $\chi <0$ (they are the quotient singularities, so are log-terminal) and $\chi =0$ ( a well-known list of the log-canonicals).  
%\end{remark}
\begin{corollary}  \label{co:g} Suppose $X$ is a $\Q$-Gorenstein weighted homogeneous surface singularity (e.g., $g=0$), not  a cyclic quotient.  Then the graded discrepancy is
 $$\alpha(X)=-1-(\chi/e).$$
 %By a previous remark, the largest coefficient of $Z$ is $\chi/e$.
\begin{proof}  The Seifert partial resolution is obtained from the minimal good resolution by blowing-down the $t$ strings.  The parenthetical assertion in the Corollary, that a weighted homogeneous singularity with rational central curve is $\Q$-Gorenstein, follows for instance from the main theorem of \cite{n}; this implies the stronger result that such an $X$ is a quotient of a complete intersection.
\end{proof}
\end{corollary}

Next, it is natural to find the numerical (i.e., topological) order of the class of $K_X$ (or equivalently $K_{\tilde{X}}+E$) in the discriminant group $\Delta=\mathbb E^*/\mathbb E$ of the singularity.  A multiple $rZ$ is integral if and only if its image in $\Delta$ dots to $0$ with every element in the discriminant group, or  with any set of generators.  But, the duals of any $t-1$ of the end curves are generators, and $Z$ dotted with an end curve dual is equal to the (negative of the) coefficient of the corresponding end curve in $Z$.  At the end of the $k$th string, this is 
                                $$(\beta -q_k')/n_k.$$
  We summarize in
  
  \begin{proposition}  Consider the graph of a weighted homogeneous surface singularity as above. Then
  \begin{enumerate}
  \item The coefficient of $C$ in $-K$ is $1+(\chi/e).$
  \item All other coefficients of $-K$ are strictly smaller than $1+(\chi/e).$
  %\item One has $$-(K+E)^2=\chi^2/e\ +\sum_{k=1}^t q_k'/n_k.$$
%  \item One has on the resolution of the singularity $$c_1^2+c_2=(K+E)^2-\sum_i (d_i-3)-6g+3,$$ the sum over every curve in $E$.
  \item The order of the class of $K$ in the discriminant group is the least common denominator for any $t-1$ of the terms $$(\chi/e -q_k')/n_k.$$ 
  \end{enumerate}
  \end{proposition}
  
 The second assertion must be checked separately for quotient singularities.
 
 Combining Corollary \ref{co:g}  with Corollary \ref{co:f}, one sees there is a special role played by singularities with $\chi/e < 1.$  Observe that $q_i=1$ means the corresponding string has length $1$.  The following result is easily verified from the definitions of $\chi$ and $e$.
 
 \begin{lemma} Consider a weighted homogeneous surface singularity with rational central curve, of self-intersection $-d$, and Seifert invariants $(n_i,q_i), 1\leq i\leq t.$  Then $\chi/e <1$ in each of the following cases:
 \begin{enumerate}
 \item $t=3$ and either

\begin{enumerate}
 \item $d\geq 4$
 \item $d=3, q_1=1$
 \item $d=2, q_1=q_2=1.$ 
 \end{enumerate}
 \item $t=4,\  d \geq 3$, and $q_1=q_2=q_3=1.$
 \end{enumerate}
 \end{lemma}

\section{The main theorem} 
Let $(X,0)$ be the germ of a normal surface singularity, topologically the cone over its neighborhood boundary $L$, its link.  $L$ is a compact three-manifold, and can be viewed as the boundary of a tubular neighborhood of the exceptional fibre in a resolution of $X$.  If now $f:(\mathcal X,0)\rightarrow (\C,0)$ is a smoothing, then the general fibre $M$ of $f$, the \emph{Milnor fibre}, is a compact four-manifold with boundary $L$.  $M$ has the homotopy type of a complex of dimension $2$, and its first betti number is $0$.  If $N$ denotes the link of $\mathcal X$, then $N$ is a compact five-manifold containing $L$, whose complement fibres over the circle via $f/|f|$, with general fibre $M$.  (This is called an \emph{open-book decomposition of $N$}).  In particular, $N$ can be constructed from $M$ by adjoining cells of dimension $>2$, so that they have the same fundamental group, and  $H_2(M)\rightarrow H_2(N)$ is surjective.  For further details, see \cite{lw},  Lemma 5.1.  

Before stating and proving the main theorem, we start with a general

\begin{proposition}  Let $f:(\mathcal X,0)\rightarrow (\C,0)$ be a smoothing of a normal surface singularity $(X,0)$ so that the Milnor fibre $M$ is a rational homology disk (i.e, $b_2(M)=0$).  Then
\begin{enumerate}
\item $(X,0)$ is a rational surface singularity
\item $(\mathcal X,0)$ is a rational three-fold singularity
\item The link $N$ of $\mathcal X$ is a rational homology sphere
%\item $H_1(N)\cong H_1(M)\cong H_1(\partial M)/I$, where $I$ is a self-isotropic subgroup for the discriminant pairing on the first homology of the link $\partial M$ of  $X$.
\item Pic $(\mathcal X-\{0\})$ is finite, hence $K_{\mathcal X}$ is $\Q$-Cartier.
\end{enumerate}
\begin{proof}  That $(X,0)$ is rational follows because of the Durfee-Steenbrink formula $\mu_0+\mu_+=2p_g(X)$ (\cite{st}).  R. Elkik's theorem \cite{e}  then implies that $(\mathcal X,0)$ is rational.  As mentioned above, $H_1(M;\Q)\cong H_1(N;\Q)=0$ in general, and $H_2(N;\Q)=0$ because $M$ has no rational homology.  Thus, the compact five-manifold is a rational homology sphere.

Following Mumford's original arguments in \cite{mu}, let $U=\mathcal X-\{0\}$, and consider the cohomology sequence $$H^1(U,\mathcal O_U)\rightarrow H^1(U,\mathcal O^*_U)\rightarrow H^2(U,\mathbb Z).$$
The first term is $0$ because the depth of $\mathcal X$ is $3$.  The five-manifold $N$ is homotopic to $U$, so $H^2(U,\mathbb Z)\cong H^2(N,\mathbb Z)$ is finite; this is then isomorphic to the torsion in $H_1(N;\mathbb Z)$, which is in turn isomorphic to $H_1$ of the Milnor fibre.  In particular, Pic($U$) is finite.
     
\end{proof}
\end{proposition}
\begin{remark}   Artin component smoothings of a rational surface singularity are definitely not $\Q$-Gorenstein (except for the rational double points); for, the Milnor fibre, being diffeomorphic to the minimal resolution, is simply connected, so is not the quotient  of the Milnor fibre of the index one cover.  In fact, the total space $\mathcal X$ of such a smoothing does not have $K_{\mathcal X}$ $\Q$-Cartier.
\end{remark}
\begin{remark}  As mentioned in the Introduction, it does not follow immediately from the Proposition that the total space $\mathcal X$ is $\Q$-Gorenstein, since it is not clear that its index $1$ cover $\mathcal Y\rightarrow \mathcal X$ is Cohen-Macaulay.  In particular, $\mathcal Y$ could conceivably give a smoothing of a non-normal model of the index one cover of $X$. \end{remark}

 We  now  state and prove the main theorem of this paper.
 
 \begin{theorem}  Let $(X,0)$ be a weighted homogeneous surface singularity, possessing a rational homology disk smoothing.  Then
 \begin{enumerate} 
 \item The smoothing is induced by the graded smoothing $f:(\mathcal X,0)\rightarrow (\C,0)$ coming from a one-dimensional smoothing component in the base space of the semi-universal deformation of $(X,0)$.
 \item $(\mathcal X,0)$ is log-terminal, and the smoothing is $\Q$-Gorenstein.
 \end{enumerate}
 \begin{proof}
The result is known for the relevant cyclic quotients by explicit construction, so we skip this case.  The grading on $X$ extends to a grading on the base-space and total space of the semi-universal deformation.   According to Corollary 8.2 of \cite{w5}, a rational homology disk smoothing occurs over a one-dimensional smoothing component; it must also be graded.   We consider the induced smoothing $f:(\mathcal X,0)\rightarrow (\C,0)$ over the normalization of this component.  By Proposition 3.1, $K_{\mathcal X}$ is $\Q$-Cartier.  By Corollary 2.5 of \cite{SSW}, the resolution dual graph of $X$ satisfies one of the  cases (1) or (2) of Lemma 2.4; therefore, $\chi/e <1.$ By Corollary 2.2, it follows that $\alpha(X)>-2.$  Corollary 1.7 now yields the desired result.
 \end{proof}
 \end{theorem}

 \begin{remark}
 \begin{enumerate}
 \item As J. Koll\'ar has pointed out, the log-terminality result need not hold after base-change of the deformation, as the condition on the weights in Proposition \ref {pr:wts} need not hold.   For instance, one would not want to consider the deformation of the elliptic cone given by $x^3+y^3+z^3+t^{3m}=0$; while the total space is a graded ring, its weights $(m,m,m,1)$ do not induce the original weights $(1,1,1)$ when one sets $t=0$.
 \item It is important to note that among all the $X$ which admit rational homology disk smoothings, the only log-canonical examples are the cyclic quotients plus three more (which are quotients of elliptic cones).
 \item The proof did not use the precise list of resolution graphs found in \cite{bs}, but only the restrictions on the graph  found earlier in \cite{SSW}.
 \end{enumerate}
 
 \end{remark} 

 \section{When is $X$ $\Q$-Gorenstein?}

The topological data of a weighted homogeneous surface singularity $X$ is given by its graph, i.e. the data $(g,d,\{n_i/q_i\})$.  (We exclude cyclic quotients.)  The additional information of the analytic type is: the central curve $C$, if $g>0$; the isomorphism class of a conormal divisor $D$ of $C$ of the resolution, which has degree $d$; and the $t$ points $P_1,\cdots,P_t$ on $C$.  A theorem of Pinkham \cite{p2} (also due to Dolgachev) shows that the graph plus this analytic data allows one to write down the graded pieces of the ring of the singularity.  We describe this result following the approach of M. Demazure \cite{d}, who proved a more general result describing graded normal domains of any dimension.  

Consider a $\Q$-divisor on $C$, say  $F=\Sigma r_jQ_j$ , where $r_j\in \Q$ and $\ Q_j\in C$.  Two such divisors are equivalent if their difference is an integral divisor, linearly equivalent to $0$.  Define the integral divisor  $$\lfloor F \rfloor=\Sigma \lfloor r_j \rfloor Q_j,$$ and the invertible sheaf $$\mathcal O(F)\equiv \mathcal O(\lfloor F \rfloor)\subset k(C).$$  (As usual,  $\lfloor r \rfloor$  means greatest integer $\leq r$).  

Now let $X=\text{Spec} \ A$ be a weighted homogeneous surface singularity with graph $\Gamma$ and analytic data $C,D,\{P_1,\cdots,P_t\}.$  Define the $\Q$-divisor $$E=D-\Sigma(q_i/n_i)P_i.$$
The degree of $E$ is the topological invariant $e>0$.  

\begin{theorem}\cite{p2} Let $X=\text{Spec}\ A$ be a graded normal surface singularity with resolution graph $\Gamma$ and analytic invariants as above.   Then
 $$A=\oplus_{k=0}^{\infty}A_k=\oplus_{k=0}^{\infty}H^0(\mathcal O(kE))T^k\subset k(C)[T].$$
\end{theorem}

(In the notation of \cite{p2} and \cite{w4}, one has  $\lfloor kE\rfloor=D^{(k)}$.)
%Then the order of the \emph{discriminant group} is $n_1\cdots n_te.$
%Recall also that in a covering of degree $m$,  $$-P\cdot P=\chi ^2/e$$ multiplies by $m$. 

In \cite{wat}, K. Watanabe computed the graded local cohomology of $A$ (in all dimensions), hence the canonical sheaf $K_X$ of $X$.  Consider the $\Q$-divisor $$\Xi=K+\Sigma (1-1/n_i)P_i,$$
where $K$ is a canonical divisor on $C$.  The degree of $\Xi$ is the topological invariant $\chi$.  

\begin{theorem} \cite{wat},(2.8).  With the notation above, the dualizing module of $A$ is 
$$\omega_A=\oplus_{k=-\infty}^{\infty} H^0(C,\Xi +kE)T^k.$$
\end{theorem}
\begin{corollary} [See also \cite{w4}, (2.1)]  $A$ is Gorenstein iff there is a $t\in \Z$ with $$\Xi\equiv tE,$$ i.e., $tq_i\equiv1 (n_i),$ all $i$, and $\lfloor tE\rfloor \equiv K_C.$  Necessarily, $t=\chi/e.$
\end{corollary}
Watanabe's method yields a formula as well for $sK_X$, i.e., the double-dual $\omega_A^{(s)}$ of the $s^{th}$ tensor power of $\omega_A$ (this is first stated in
 \cite{wat2}, Lemma 3.2).  By \cite{wat}, Theorem 1.6, there is an isomorphism of the group of equivalence classes of certain $\Q$-divisors whose fractional part involves only the $P_i$, onto the divisor class group of $A$; this map sends $$F\mapsto \oplus_{k=-\infty}^{\infty} H^0(C,F +kE)T^k.$$
 
\begin{theorem}[Watanabe]  Notation as above, $$\omega^{(s)}_A=\oplus_{k=-\infty}^{\infty} H^0(C,s\Xi +kE)T^k.$$
\end{theorem}
\begin{corollary}$sK_X$ is Cartier iff there is a $t\in \Z$ with $$s\Xi\equiv tE.$$ 
\end{corollary}
Taking degrees, the last analytic condition implies $s\chi=te$, from which one can easily prove
\begin{corollary} $X$ is $\Q$-Gorenstein iff the class of the degree $0$ $\Q$-divisor $$ \Xi-(\chi/e)E$$ is torsion.
\end{corollary}

\begin{remark}
\begin{enumerate}
\item If $X$ is $\Q$-Gorenstein, then the order of $K_X$ in the (analytic) divisor class group is divisible by the denominator  of $\chi/e$, written as a reduced fraction (since $t/s=\chi/e$).  Compare with Proposition 2.3(3) above.
\item Corollary 4.6 gives an explicit condition on a conormal divisor $D$ so that a weighted homogeneous singularity is $\Q$-Gorenstein, once one fixes the topological data plus the isomorphism type of the curve $C$ and the points $P_i$.  Recall that Popescu-Pampu has proved \cite{pp}  that for the resolution graph of \emph{any} normal surface, with fixed analytic type of the reduced exceptional divisor, one can choose normal bundles appropriately so that there exists a $\Q$-Gorenstein singularity with the given topological plus analytic data. 
\item An alternate description of a weighted homogeneous surface singularity $A$, due to I. Dolgachev, writes $$A=\oplus_{k=0}^{\infty}H^0(U,L^{-k})^{\Gamma},$$ where $L$ is a line bundle on the universal covering space $U$ of Proj $A$, and  $\Gamma\subset \text{Aut}\ U$ is a discrete subgroup.  From this point of view, the Gorenstein condition is found by Dolgachev in \cite{do}, and the $\Q$-Gorenstein condition is due to A. Pratoussevitch \cite{pr}.
\end{enumerate}
\end{remark}

 \end{document}